\documentclass{amsart}
\usepackage[utf8]{inputenc}
\usepackage{amssymb}
\usepackage{tikz-cd} 
\usepackage{amsfonts}
\usepackage{amsmath}
\usepackage{amstext}
\usepackage{amsthm}
\usepackage{latexsym}
\usepackage[shortlabels]{enumitem}

\usepackage[utf8]{inputenc}

\usepackage{mathrsfs}  
\usepackage[mathscr]{eucal}

\usepackage{amscd}
\usepackage{eucal}
\usepackage{amsxtra}
\usepackage{amsbsy}
\usepackage{graphicx}
\usepackage{latexsym}
\usepackage{pb-diagram}
\usepackage{amsopn}
\usepackage{delarray}
\usepackage{psfrag}
\usepackage{verbatim}
\usepackage{bussproofs}
\usepackage[T1]{fontenc}
\usepackage{color,cancel}

\usepackage{appendix}
\usepackage{color}
\usepackage{epigraph}
\usepackage{lipsum}
\usepackage{hyperref} 
\hypersetup{pdfpagemode={UseOutlines},
bookmarksopen=true,bookmarksopenlevel=0,hypertexnames=false,colorlinks=true,citecolor=black,linkcolor=black,urlcolor=black,pdfstartview={FitV},unicode,breaklinks=true}
\usepackage{graphicx}
\usepackage{enumitem}

\newtheorem{thm}{Theorem}[section]
\newtheorem{propo}[thm]{Proposition}
\newtheorem{coro}[thm]{Corollary}
\newtheorem{lem}[thm]{Lemma}
\newtheorem{que}[thm]{Question}
\newtheorem*{que*}{Question}
\newtheorem{problem}[thm]{Problem}
\newtheorem*{problem*}{Problem}

\theoremstyle{definition}
\newtheorem{dfn}[thm]{Definition}

\theoremstyle{remark}

\newtheorem{nota}[thm]{Notation}

\newcommand{\conc}{^\smallfrown}	
\newcommand{\rest}{\upharpoonright} 

\def\finb{{\textsf{FIN}^\cB}}

 \def\cB{{\mathcal{B}}} \def\cC{{\mathcal{C}}} \def\cD{{\mathcal{D}}}  \def\cF{{\mathcal{F}}} \def\cG{{\mathcal{G}}}  \def\cI{{\mathcal{I}}}     \def\cN{{\mathcal{N}}}  \def\cP{{\mathcal{P}}}      \def\cU{{\mathcal{U}}}    

                 \def\bR{{\mathbb{R}}}

\def\b{{\mathfrak{b}}}
\def\c{{\mathfrak{c}}}

\def\M{{\mathfrak{M}}}

\newcommand{\w}{\omega}

\title[High dimensional countable compactness]{High dimensional countable compactness and ultrafilters}

\author{C. Corral}
\address{York University, 4700 Keele St, Toronto, ON, Canada, M3J 1P3.}
\email{cicorral@yorku.ca}
\thanks{The first author acknowledges support from York University and the Fields Institute.}

\author[Memarpanahi]{P. Memarpanahi}
\address{University of Toronto Scarborough 1265 Military Trail, Toronto, Ontario M1C 1A4, Canada}
\email{pourya.memarpanahi@utoronto.ca}

\author{P. Szeptycki}
\address{Department of Mathematics and Statistics, York University, Toronto, Ontario M3J 1P3, Canada}
\email{szeptyck@yorku.ca}
\thanks{The third author acknowledges support from NSERC grant 503970.}

\date{}



\begin{document}

\maketitle

\begin{abstract}
We define several notions of a limit point on sequences with domain a barrier in $[\omega]^{<\omega}$ focusing on the two dimensional case $[\omega]^2$. 
By exploring some natural candidates, we show that countable compactness has a number of generalizations in terms of limits of high dimensional sequences and define a particular notion of $\alpha$-countable compactness for $\alpha\leq\w_1$. 
We then focus on dimension 2 and compare 2-countable compactness with notions previously studied in the literature. 
We present a number of counterexamples showing that these classes are different. In particular assuming the existence of a Ramsey ultrafilter, a subspace of $\beta\omega$ which is doubly countably compact whose square is not countably compact, answering a question of T. Banakh, S. Dimitrova and O. Gutik \cite{Banakhsemigroups}. The analysis of this construction leads to some possibly new types of ultrafilters related to discrete, P-points and Ramsey ultrafilters. 
\end{abstract}


\section{Introduction and notation}

A sequence in a topological space $X$ is a function $f:\w\to X$. Our main object of study is the generalization of a sequence to higher dimensions. For example, A double sequences is a function $f:[\omega]^2\rightarrow X$ and an $n$-dimensional sequence is a function with domain $[\omega]^n$. And more generally, for a barrier $\cB\subseteq [\omega]^{<\omega}$ (see below for the relevant definitions), we will refer to $f:\cB\to X$ as a $\cB$-sequence. 

Given a barrier $\cB$, $M\in[\w]^\w$ and $f:\cB\to X$, we say that $f$ converges to $x\in X$ if for every $U\in\cN(x)$ there exists $n\in\w$ such that $f(s)\in U$ for every $s\in\cB|(M\setminus n)$. Following \cite{corral2023infinite}, we say that a space $X$ is $\cB$-sequentially compact if for every function $f:\cB\to X$ there exists $M\in[\w]^\w$ such that $f\rest(\cB|M)$ converges to some $x\in X$. Moreover, $X$ is $\alpha$-sequentially compact if it is $\cB$-sequentially compact for every barrier $\cB$ of rank $\alpha$ and it is $\w_1$-sequentially compact if it is $\alpha$-sequentially compact for every $\alpha<\w_1$. The case $\alpha<\omega$ was first defined and studied by W. Kubi\'s and P. Szeptycki \cite{kubistopologicalramsey}, the case $\alpha=2$ was considered even earlier by M. Boja\'{n}czyk, E. Kopczy\'{n}ski and S. Toru\'{n}czyk \cite{ramseymetric}, and the infinite dimensional case was introduced and studied in 
\cite{corral2023infinite}.

The notion of $\cB$-sequential compactness has two main ingredients: On the one hand, its compactness nature ensures the existence of limit points for countable sets of the form $f[\cB]$ where $f:\cB\to X$. On the other hand, dimension is captured by the structure of the barrier $\cB$ and the mode of convergence attached to $\cB$-sequences. Even if we drop one of the ingredients, the resulting notions seem to be very natural properties for a topological space. 
Indeed, by dropping compactness from the notion of 2-sequential compactness, one obtains the Ramsey property considered by H. Knaust (\cite{knaustangelic},\cite{knaustarray}) and if instead one focuses on the compactness side, there are a number of natural properties generalizing well known notions. In particular, we will be interested in high dimensional versions of countable compactness. 
We will motivate and define the notion of $\cB$-countable compactness and see how it relates with other notions of high dimensional countable compactness considered in the past, like the notion of \emph{doubly countably compact} introduced by T. Banakh, S. Dimitrova and O. Gutik \cite{Banakhsemigroups}.

Various counterexamples delineating these properties are presented, some of which depend on the existence of a Ramsey ultrafilter in $\omega^*$. An analysis of these constructions leads us to the formulation of a new types of special points in $\omega^*$ weaker than Ramsey and related to Baumgartner's notion of a discrete ultrafilter (see section 4 below).

\vspace{5mm}

Our notation and terminology are fairly standard.  For a set $A$ and a cardinal $\kappa$, we denote by $[A]^\kappa$ the family of subsets of $A$ of size $\kappa$ and $[A]^{<\kappa}=\bigcup_{\lambda<\kappa}[A]^\lambda$. An ideal $\cI\subseteq \cP(X)$, is a nonempty family closed under subsets and finite unions. We will always assume that all finite subsets of $X$ are members of $\cI$ and that $\cI$ is proper, i.e., $X\notin\cI$. Similarly, $\cF\subseteq \cP(X)$ is a filter if $\cF^*=\{X\setminus A:A\in\cF\}$ is an ideal. Following our previous convention, every filter contains the cofinite subsets of $X$, however, it will be useful to overlook this assumption sometimes, so in the particular case of $X=\w$, we say that $\cF$ is a free ultrafilter or a non-principal ultrafilter if it contains all cofinite subsets. A filter $\cF$ is an ultrafilter if either $A$ or $X\setminus A$ is in $\cF$ for every $A\subseteq X$. Given $\cG\subseteq\cP(X)$, the filter generated by $\cG$, denoted by $\langle\cG\rangle$, is the family of all subsets of $X$ that contain the intersection of finitely many elements from $\cG$.

For a function $f$ and a subset $M$ of its domain, we denote by $f''M$ or $f[M]$ its image. We mainly use the notation $f''M$ when $f$ is defined on some power of a superset of $M$, for instance, if $f:[\w]^2\to X$ and $M\subseteq\w$, we denote by $f''[M]^2$ the image of $[M]^2$ under $f$.
For $A,B\subseteq\w$, we say that $A$ is almost contained in $B$ and denote it by $A\subseteq^*B$, if $A\setminus B$ is finite. Given $\cG\subseteq[\w]^\w$, a set $P\in[\w]^\w$ is a pseudointersection for $\cG$ if $P\subseteq^*G$ for every $G\in\cG$.

A tree $T$ is a partially ordered set such that $prec(t)$ is well ordered for every $t\in T$, where $prec(t)=\{s\in T:s\leq t\}$. We say that $T$ is \emph{well pruned} if for every $t\in T$ there is $s\in T$ such that $t\leq s$ and $t\neq s$. In the case where $T\subseteq\w^{<\w}$, we fix the notation $succ(t)=\{n\in\w:t\conc n\in T\}$, for the set of successors of $t$ in $T$.

For a topological space $X$ and $x\in X$, we denote by $\cN(x)$ the set of open neighborhoods of $x$. We usually identify sequences with their range. A point $x\in X$ is a limit point of $A$ if every $U\in\cN(x)$ intersects $A$ in an infinite set and it is an accumulation point for $A$ if $U\setminus\{x\}\cap A\neq\emptyset$ for every $U\in\cN(x)$.
Given a free ultrafilter $p$ on $\w$ and a sequence $\{x_n:n\in\w\}$, we say that $x$ is the $p$-limit of the sequence if $\{n\in\w:x_n\in U\}\in p$ for every $U\in\cN(x)$. It is clear that $p$-limits are unique in Hausdorff spaces (and in fact, all spaces considered here are Tychonoff). A space $X$ is $p$-compact if every sequence in $X$ has a $p$-limit. 
It is well known that $p$-compactness sits between compactness and countable compactness for every free ultrafilter $p$ (see \cite{bernsteinnewkind}). 

The notion of filter convergence was introduced by Bourbaki \cite{bourbakifilters} (see also \cite{engelking}). For a filter $\cF\subseteq \cP(X)$ and $x\in X$, we say that $\cF$ clusters at $x$ if $x\in \overline{F}$ for every $F\in\cF$ and $\cF$ converges to $x$ if $\cN(x)\subseteq \cF$. If $\cF$ is a filter base, then we use the same terminology if the filter generated by $\cF$ converges or clusters. A topological space $X$ is \emph{bisequential at} $x$ if for every filter $\cF$ that clusters at $x$ there is a countable family $\cG=\{F_n:n\in\w\}\subseteq\cF$ such that $\langle\cG\rangle$ converges to $x$. The whole space $X$ is \emph{bisequential} if it is bisequential at every $x\in X$.

Given a filter $\cF$ on a countable set $N$ and $f:N\to X$, we will also say that $\cF$ converges to (clusters at) $x\in X$ if the filter base $f(\cF)$ converges to (clusters at) $x$.

The Stone-\v{C}ech compactification of $\w$ will be denoted by $\beta\w$ and can be identified with the space of all ultrafilters on $\w$. On $\beta\w$, the closure of every countable discrete subset is homeomorphic to $\beta\w$, and $x\in\overline{\{x_n:n\in\w\}}$, where $\{x_n:n\in\w\}$ is discrete, if and only if $x$ is the $p$-limit of that sequence for some ultrafilter $p$ (if $p$ is not free then $x=x_n$ for some $n$).

\vspace{3mm}

For $\cB\subseteq[\omega]^{<\omega}$ and $M\in[\omega]^\omega$ let $\cB|M=\{s\in\cB:s\subseteq M\}$.
The relation $s\sqsubseteq t$ means that $s$ is an initial segment of $t$ and $s\sqsubset t$ means $s\sqsubseteq t$ and $s\neq t$. We say that $\cB\subseteq[M]^{<\w}$ is a \emph{barrier}, where $M\in[\w]^\w$ if it is an antichain with respect to $\sqsubseteq$ and every infinite subset of $M$ contains an initial segment in $\cB$. We assume that a barrier $\cB$ consists always of nonempty sets and it is nontrivial if $\cB\neq[\w]^1$.

Given a barrier $\cB$, let 
$$T(\cB)=\{s\in[\omega]^{<\omega}:\exists t\in\cB\ (s\sqsubseteq t)\}$$ 
and $\rho_{T(\cB)}:T(\cB)\to\omega_1$ be given by
$$\rho_{T(\cB)}(s)=\sup\{\rho_{T(\cB)}(t)+1:t\in T(\cB)\land s\sqsubset t\},$$
where $\sup\emptyset=0$. We will omit the subindex $T(\cB)$ when there is no risk of confusion. The \emph{rank of} $\cB$ is defined as $\rho(\cB)=\rho_{T(\cB)}(\emptyset)$.

For a given barrier $\cB$ and $s\in[\w]^{<\w}$, we let $\cB(s)=\{b\in\cB:s\sqsubseteq b\}$ and $\cB[s]=\{b\setminus s:b\in\cB(s)\}$. We will identify $[\w]^{<\w}$ with $\w^{<\w}$ by taking the increasing enumeration of a finite set, mainly when dealing with barriers $\cB$ and their associated trees $T(\cB)$.

We also define an ideal naturally associated to a barrier $\cB$ by recursion on the rank of $\cB$. If $\cB=[\w]^1$, then $\textsf{FIN}^\cB=\textsf{FIN}$ is the ideal of finite subsets of $\{\{n\}:n\in\w\}$. For an arbitrary barrier $\cB$, $A\in\textsf{FIN}^\cB$ if and only if $A[n]\cap\cB[\{n\}]\in\textsf{FIN}^{\cB[\{n\}]}$ for all but finitely many $n\in\w$, where $A[n]=\{s\setminus\{n\}:s\in A\land \min(s)=n\}$.
When $\cB=[\w]^n$, this ideal, usually denoted by $\textsf{FIN}^n$ has been widely studied. 

\section{Limit points and countable compactness}

In this section we will deal with different kinds of limit points for high dimensional sequences. In dimension 1, many notions like limit points, cluster points and accumulation points do coincide for Hausdorff spaces. However, these notions though equivalent, differ formally, and they give rise to different kinds of limit points for higher dimensions. 
We will discuss how most of these properties are related (and sometimes coincide) with classical notions and subsequently we will use them to define high dimensional compactness notions. In particular, we will motivate our choice for the definition of high dimensional countable compactness.\\

For a free ultrafilter $p$ on $\w$ and a well-pruned tree $T\subseteq \w^{<\w}$, we say that $T$ is a $p$\textit{-branching tree} if $succ_T(t)\in p$ for every $t\in T$.
Moreover, if $\cB$ is a barrier on $\w$, we define an ultrafilter $$p^\cB=\{T\cap \cB:T\textnormal{ is a $p$-branching tree}\}\subseteq\cP(\cB).$$

In particular, this ultrafilter coincides with $p^k$ as defined in \cite{todorcevicramseyspaces} when $\cB=[\w]^k$.

\begin{dfn}\label{limit points}
    Let $X$ be a topological space, let $p$ be a free ultrafilter on $\w$, let $\cB$ be a barrier, $f:\cB\to X$ and $x\in X$. We say that $x$ is:
\begin{itemize}
    \item\label{limititem1} the $\textsf{FIN}^\cB$\emph{-limit point of $f$} if the set $\{s\in \cB:f(s)\notin U\}\in\finb$ for every $U\in\cN(x)$.

    \item\label{limititem2} the $p^\cB$\emph{-limit of} $f$ if $f^{-1}(U)\in p^\cB$ for every $U\in\cN(x)$.

    \item\label{limititem3} A $\cB$-accumulation point for $f$ if for every $U\in\cN(x)$ there exists $M\in[\w]^\w$ such that $f[\cB|M]\subseteq U$.
    
    \item\label{limititem4} A $\cB$\emph{-limit point} of $f$ if $f^{-1}(U)\notin\finb$ for every $U\in\cN(x)$.\\
\end{itemize}
\end{dfn}

The first and fourth notions can be restated in terms of filter convergence. 
Let $\cF$ be the dual filter to $\finb$. Then $x$ is the $\finb$-limit of $f:\cB\to X$ if and only if the filter $\cF$ converges to $x$, and $x$ is a $\cB$-limit point for $f$ if and only if $\cF$ clusters at $x$. In particular, if $\cB=[\w]^1$, we can identify $f$ with a sequence $\{x_n:n\in\w\}$, thus $x$ is the $\finb$-limit of $f$ if and only if $\{x_n:n\in\w\}$ converges to $x$ and $x$ is a $\cB$-limit point for $f$ if and only if it is a limit point for the sequence $\{x_n:n\in\w\}$.

It is also worth mentioning that $\finb$-limits coincide with a type of convergence considered by M. Balcerzak and K. Dems in \cite{Balcerzaktypesofconvergence} for double (and $n$-dimensional) sequences of real numbers.

\begin{propo}\label{limitsimplications}
Given a $\cB$-sequence $f:\cB\to X$, $x\in X$ and a free ultrafilter $p$ we have that:
$$f\textnormal{ converges to }x\implies x\textnormal{ is the }\finb\textnormal{-limit of }f\implies x\textnormal{ is the }p^\cB\textnormal{-limit of }f$$
$$\implies x\textnormal{ is an accumulation point for }f\implies x\textnormal{ is the }\cB\textnormal{-limit of }f$$
\end{propo}

\begin{proof}
    Let $U\in\cN(x)$. If $f$ converges to $x$ then $\cB|(\w\setminus k)\subseteq f^{-1}(U)$ for some  $k\in\w$, which clearly shows that the complement of $f^{-1}(U)$ is an element of $\finb$.\\
    In turn, this implies that $f^{-1}(U)\in p^\cB$ as the dual filter of $\finb$ is contained in $p^\cB$ for every free ultrafilter $p$.\\ 
    Assume now that $f^{-1}(U)\in p^\cB$ and define 
    $$T(U)=\{s\in T(\cB):\exists b\in f^{-1}(U):s\sqsubseteq b\}.$$ 
    We can define a set $M=\{k_n:n\in\w\}$ by choosing $k_n$ in such a way that if $s\subseteq\{k_i:i<n\}\cap T(\cB)\setminus\cB$ then $s\conc k_n\in T(U)$. This is possible as $succ_{T(U)}(s)\in p$ for every $s\in T(U)$. Therefore we have that $f[\cB|M]\subseteq U$.\\
    Finally, it is clear that any set of the form $\cB|M$ is not an element of $\finb$ and so any $\cB$-accumulation point for $f$ is also a limit point of $f$.
\end{proof}

As any two elements in $p^\cB$ have non-empty intersection, $p^\cB$ limit points are unique in Hausdorff spaces (and so are $\finb$-limit points in consequence). On the other hand, if $M_0$ and $M_1$ are infinite disjoint subsets of $\w$ and $x_0\neq x_1$, we can define a $\cB$-sequence $f$ such that $f\rest(\cB|M_i)$ is constant with value $x_i$, and then $x_i$ is a $\cB$-accumulation point for every $i\in\{0,1\}$, that is, $\cB$-accumulation points (and then $\cB$-limit points) are not unique. This justifies our use of the definite and indefinite articles ``the'' and ``a'' in definition \ref{limit points}.

\vspace{3mm}

It is a basic fact in topology, that a space $X$ is countably compact if any of the following equivalent conditions holds:
\begin{enumerate}
    \item Every countable cover of $X$ has a finite subcover,
    \item every one-to-one sequence has an accumulation point and
    \item every sequence has a $p$-limit for some free ultrafilter $p$ on $\w$.
\end{enumerate}

We now introduce the compactness properties associated to the limit points defined above. Our choice for the high dimensional version of countable compactness will be justified later in the paper.

\begin{dfn}Given a space $X$ and a barrier $\cB$ we say that:
    \begin{enumerate}
            \item $X$ is $\finb$\emph{-compact} if for every function $f:\cB\to X$ there is an infinite set $M$ such that $f\upharpoonright (\cB|M)$ has a $\finb$-limit point.\footnote{Formally, we should say ``a $\textsf{FIN}^{(\cB|M)}$-limit point'', but we will soon get rid of this property and we see no point on being so pedantic.}

            \item\label{itemcountablycompact} $X$ is $\cB$-\emph{countably compact} if for every $f:\cB\to X$ there exists a free ultrafilter $p$ on $\w$ such that the $p^\cB$-limit of $f$ exists.
            
            \item $X$ is \emph{weakly} $\cB$\emph{-countably compact} if every function $f:\cB\to X$ has a $\cB$-accumulation point.

            \item\label{itemweaklycountablycompact} $X$ is $\cB$-\emph{limit point compact} if every $f:\cB\to X$ has a limit point.
    \end{enumerate}
\end{dfn}

Note that for item (\ref{itemcountablycompact}), by fixing a bijection $\varphi:\cB\to\w$, we can identify $\cB$-sequences with classical one dimensional sequences, and the ultrafilter $p^\cB$ corresponds to an ultrafilter $q=\varphi(p^\cB)$. Hence $\cB$-countable compactness is equivalent to the assertion that every sequence has a $q$-limit for some free ultrafilter $q$ on $\w$ which is of the form $\varphi(p^\cB)$. 
Let's denote by $\cU_\varphi(\cB)$ this family of ultrafilters (which of course depends on $\varphi$). 
It is easy to see that $\cU_\varphi(\cB)$ is a proper subset of $\w^*$, as any ultrafilter containing a set of the form $\cB(s)$ for some $s\in T(\cB)$ is not in $\cU_\varphi(\cB)$. 
Consequently, $\cB$-countable compactness is a priori, stronger than countable compactness. We will see that this is indeed the case in Section \ref{2dimensionalsection}.
The previous discussion suggests the notion of $D$-compactness for $D\subseteq\w^*$.

\begin{dfn}
    Let $D\subseteq\w^*$. A space $X$ is $D$\textit{-compact} if every sequence in $X$ has a $p$-limit for some ultrafilter $p\in D$.
\end{dfn}

A similar notion for pseudocompactness was introduced by S. Garc\'ia-Ferreira and Y. Ortiz-Castillo in \cite{garcia2018pseudocompactness}. It is now clear from definitions that being $\cB$-countably compact becomes the natural property of being $D$-compact for $D=\cU_\varphi(\cB)$.

One of the motivations for our choice of the high dimensional version of countable compactness is the following: We know every countably compact Fr\'echet\footnote{Recall that a space is Fr\'echet if every point $x$ is in the closure of a set $A$ if and only if there is a sequence in $A$ converging to $x$} space is sequentially compact. In the same spirit, it was shown in \cite{corral2023infinite} that compact bisequential spaces are $\cB$-sequentially compact for every barrier $\cB$. 
Bisequentiality is a generalization of the Fr\'echet property in the language of filters (actually every bisequential space is Fr\'echet), and the consequence of compactness used there, is that every $\cB$-sequence has a $p^\cB$-limit point for some free ultrafilter $p$. Hence we can restate that theorem with our current terminology as follows.

\begin{thm}\cite{corral2023infinite}\label{bcountablybisequentialisbseq}
    Every $\cB$-countably compact bisequential space is $\cB$-sequentially compact.\qed
\end{thm}

We first show that the class of $\finb$-compact spaces coincides with the class of sequentially compact topological spaces.

\begin{thm}\label{2.6}
    The following are equivalent:
    \begin{enumerate}[1)]
        \item $X$ is $\finb$-compact for every barrier $\cB$,
        \item $X$ is $\finb$-compact for some barrier $\cB$,
        \item $X$ is sequentially compact.
    \end{enumerate}
\end{thm}

\begin{proof}
    $1)\rightarrow2)$ is trivial, so let us start with 2) $\rightarrow$ 3): If $X$ is $\finb$-compact and $\{x_n:n\in\w\}\subseteq X$, we can define $f:\cB\to X$ in such a way that $f(s)=x_n$ whenever $n=\min(s)$. 
    If $M\in[\w]^\w$ and $x\in X$ are such that $x$ is the $\finb$-limit of $f\rest(\cB|M)$, then $\{x_n:n\in M\}$ converges to $x$.

    Finally, for $3)\rightarrow1)$, let $X$ be sequentially compact. We shall show that $X$ is $\finb$-compact by induction on the rank of the barrier $\cB$.\\
    Assume that $X$ is $\textsf{FIN}^\cC$-compact for every barrier $\cC$ such that $\rho(\cC)<\rho(\cB)$.
    Let $f:\cB\to X$ and define $f_n:\cB[\{n\}]\to X$ by $f_n(s)=f(\{n\}\cup s)$ for every $s\in\cB[\{n\}]$.
    As $X$ is $\textsf{FIN}^{\cB[\{n\}]}$-compact for every $n\in\w$, we can find a decreasing sequence $\{M_n:n\in\w\}\subseteq[\w]^\w$ and a sequence $\{x_n:n\in\w\}\subseteq X$ such that $x_n$ is the $\textsf{FIN}^{(\cB[\{n\}]|M_n)}$-limit of $f_n\rest(\cB[\{n\}]|M_n)$.\\
    Take a pseudointersection $M'$ for $\{M_n:n\in\w\}$. As $X$ is sequentially compact, there are $M\in[M']^\w$ and $x\in X$ such that $\{x_n:n\in M\}$ converges to $x$. It is now easy to see that $x$ is the $\finb$-limit of $f\rest(\cB|M)$.
\end{proof}

We can now relate sequential compactness and $p$-compactness with $\cB$-countably compact spaces.

\begin{coro}\label{seqcompactisbcountably}
    Let $X$ be a sequentially compact space, then it is $\cB$-countably compact for every barrier $\cB$.
\end{coro}

\begin{proof} Let $f:\cB\rightarrow X$ be given. By Theorem \ref{2.6} fix $x\in X$ and $M$ witnessing $\finb$-compactness. Letting $p$ be any ultrafilter containin $M$ it follows that $x$ is a $p^\cB$ limit of $f$ as required.   
\end{proof}

In a similar way, we can show that countably compact spaces can be characterized in terms of $\cB$-limit point compactness. As the proof is similar to the proof of Theorem \ref{2.6}, using countable compactness instead of sequential compactness, we omit it.

\begin{thm}\label{countablycompactBlimit}
    The following are equivalent:
    \begin{enumerate}[1)]
        \item $X$ is $\cB$-limit point compact for every barrier $\cB$,
        \item $X$ is $\cB$-limit point compact for some barrier $\cB$,
        \item $X$ is countably compact.\qed
    \end{enumerate}
\end{thm}

\begin{thm}
   Suppose $X$ is $p$-compact. Then $X$ is $p^\cB$-compact for every barrier $\cB$. In particular $X$ is $\cB$-countably compact for every barrier $\cB$.
\end{thm}

\begin{proof}
    Let $X$ be $p$-compact. We will show that $X$ is $p^\cB$-compact for every barrier $\cB$ by induction on the rank of $\cB$. As $X$ being $p$-compact is equivalent to $X$ being $[\w]^1$-compact, we have the base step for free.

    Let $f:\cB\to X$ where $\alpha=\rho(\cB)>1$ and assume that $X$ is $p^\cC$-compact for every barrier $\cC$ of rank less that $\alpha$. For each $n\in\w$, define $f_n:\cB[\{n\}]\to X$ by $f_n(s)=f(\{n\}\cup s)$. By induction we can find $x_n\in X$ that is the $p^{\cB[\{n\}]}$-limit of $f_n$. Since $X$ is $p$-compact we can also find $x=p-\lim\{x_n:n\in\w\}$. It follows by a similar argument to the one in the proof of Theorem \ref{seqcompactisbcountably} that $x$ is also the $p^\cB$-limit of $f$.
\end{proof}

From the previous two propositions and Theorem \ref{bcountablybisequentialisbseq} we get the following corollary.

\begin{coro}
    Every bisequential space that is either, sequentially compact or $p$-compact, is also $\w_1$-sequentially compact.
\end{coro}

We now give connections between the remaining properties defined so far.

\begin{propo}
    $\cB$-countably compact spaces are weakly $\cB$-countably compact and weakly $\cB$-countably compact spaces are are countably compact. 
\end{propo}

\begin{proof}
    Let $X$ be a $\cB$-countably compact space and $f:\cB\to X$. Let $x\in X$ be the $p^\cB$-limit of $f$. It follows from Proposition \ref{limitsimplications} that $x$ is also an accumulation point for $f$, hence $X$ is weakly $\cB$-countably compact.

    Similarly, as any accumulation point for $f$ is also a limit point, we get that weakly $\cB$-countably compact spaces are $\cB$-limit point compact and therefore countably compact by Theorem \ref{countablycompactBlimit}.
\end{proof}

It was proved in \cite{corral2023infinite} that every barrier $\cB$ on $\w$ of finite rank $n$, is of the form $[M]^n$ on some cofinite set $M\subseteq\w$. For this reason, we can restrict our attention to barriers of the form $[\w]^n$ in the finite case. We will say that $x$ is the $p^n$-limit of $f:[\w]^n\to X$ instead of $x$ being the $p^{[\w]^n}$-limit and we will say that $X$ is $n$-countably compact if $X$ is $[\w]^n$-countably compact. In general we have the following definition. 

\begin{dfn}\label{alfacountably}
    A space $X$ is (weakly) $\alpha$\emph{-countably compact} if it is (weakly) $\cB$-countably compact for every barrier $\cB$ of rank $\alpha$.
    Moreover, a space will be said to be $\w_1$\emph{-countably compact} if it is $\alpha$-countably compact for every $\alpha<\w_1$.
\end{dfn}

{\begin{figure}[t]\label{Highdimensionalcountablecompactness}
    \centering
    \includegraphics[height=10cm]{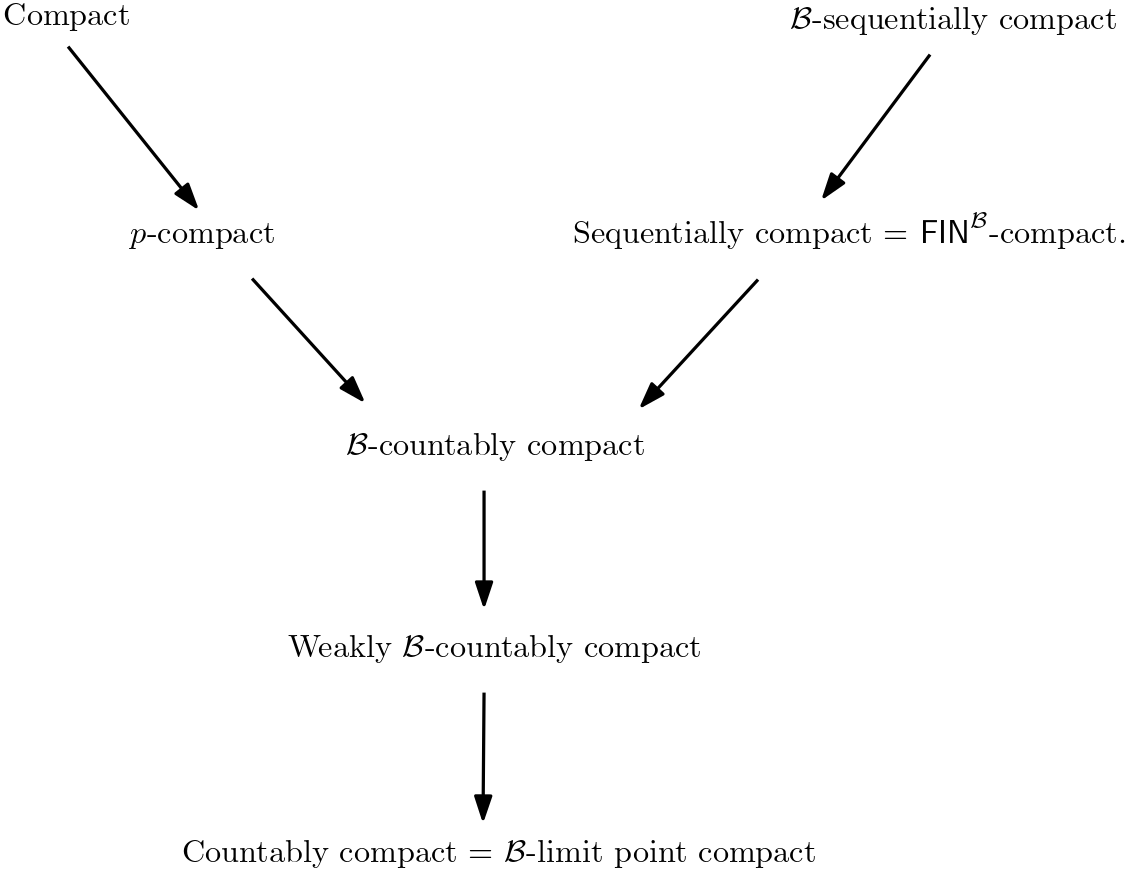}
    \caption{High dimensional countable compactness properties ($\alpha\in\w_1$).}
\end{figure}}

Figure 1 shows the relations we have proved among these high dimensional versions of countable compactness. We shall show that the arrows in the diagram do not reverse. Examples of a sequentially compact non compact space and compact non sequentially compact spaces are well known: $\beta\w$ is compact but has no nontrivial convergent sequences and $\w_1$ is even $\w_1$-sequentially compact but fails to be compact.

One of the most interesting questions in General Topology is Scarborough-Stone problem: Is the product of sequentially compact spaces countably compact?

This problem was first posed by C. T. Scarborough and A. H. Stone in \cite{ScarboroughStoneProductofnearlycompact} and listed as problem 13 in the compendium of the twenty most important problems in Set Theoretic Topology in \cite{hruvsaktwentyproblems}. The problem is known to be independent in the class of perfectly normal spaces (see \cite{vaughanproducts} and \cite{nyikoshereditary}) and it has a negative answer for Hausdorff spaces (see \cite{Nyikosscarboroughforhausdorff}). 
Hence the portion of this question that remains unsolved asks whether the product of Tychonoff sequentially compact spaces are countably compact. 
Having at hand strengthening of both sequential compactness and countable compactness, we could ask for any $\alpha,\beta\in[1,\w_1]$ whether the product of a family of $\alpha$-sequentially compact spaces is $\beta$-countably compact.

Assuming $\b=\c$, there is a family of first countable sequentially compact spaces (hence $\w_1$-sequentially compact by a result in \cite{corral2023infinite}) whose product is not countably compact. However, perhaps a weak counterexample to Scarborough Stone may exist in ZFC:

\begin{que} Is there a family of sequentially compact spaces whose product is not $\omega_1$-countably compact?
\end{que}

Figure 1 summarizes the ZFC implications between these compactness properties. Where compact implies $p$-compact for every $p$ and being $p$-compact for some $p$ implies $\cB$-countably compactness for every $\cB$.
None of the arrows reverses even consistently, with the following possible exception: We will see an example of a countably compact not $2$-countably compact but we do not know whether it is weakly $2$-countably compact or not. Also the counterexample for Scarborough-Stone mentioned above is a consistent example of an $\w_1$-sequentially compact space that is not $p$-compact for any $p$, which of course is open in {\sf ZFC}.



\section{Idempotents in semigroups and double sequences}\label{2dimensionalsection}

Most of the applications of the high dimensional versions of sequential compactness introduced in \cite{corral2023infinite} follow from the 2-dimensional case. Indeed, one of the first motivations for the introduction of n-sequential compactness in \cite{kubistopologicalramsey} is a theorem due to M. Boja\'nczyk, E. Kopczy\'nski and S. Toru\'nczyk, that asserts that compact metric semigroups contain idempotents by showing that those semigroups are 2-sequentially compact \cite{ramseymetric}. A precise characterization of the existence of idempotents in topological semigroups, in terms of a 2 dimensional version of compactness is due to T. Banakh, S. Dimitrova and O. Gutik:

\begin{dfn}\cite{Banakhsemigroups}
    Let $p$ be a free ultrafilter on $\w$. We say that $x\in X$ is the \emph{double $p$-limit} of the double sequence $\{x_{n,m}:n<m<\w\}\subseteq X$ if there are $A\in p$ and a sequence $\{x_n:n\in A\}\subseteq X$ such that $x_n$ is the $p$-limit of $\{x_{n,m}:m\in\w\setminus(n+1)\}$ for every $n\in A$ and $x$ is the $p$-limit of $\{x_n:n\in A\}$. 

    A space $X$ is \emph{doubly countably compact} if every double sequence has a double $p$-limit for some $p\in\w^*$.
\end{dfn}

Notice that if $x$ is the double $p$-limit of $\{x_{n,m}:n<m<\w\}$, then $x$ is also the $p^2$-limit of the corresponding double sequence. It will be useful later to make it explicit the set of intermediate points that turn a $p^2$-limit point into a double $p$-limit, so we will say that $x$ is the double $p$-limit  of $\{x_{n,m}:n<m<\w\}$ through $\{x_n:n\in A\}$ for $A\in p$ and $\{x_n:n\in A\}\subseteq X$ or that $(x,x_n)_{n\in A\in p}=\{x\}\cup\{x_n:n\in A\}$ is the $p$\emph{-limit sequence} of $\{x_{n,m}:n<m<\w\}$ following the notation in the definition above.

We fix some notation for the rest of this section. Let $X$ be a topological semigroup and $x\in X$. We define $f_x:[\w]^2\to X$ by $f_x(\{n,m\})=x^{m-n}$ where $n<m$.

\begin{thm}\cite{Banakhsemigroups}\label{banakhstheorem}
    A topological semigroup $X$ has an idempotent if and only if the double $p$-limit of $f_x$ exists for some $x\in X$ and $p\in\w^*$.
\end{thm}

The proof of the previous theorem actually shows that the same result holds if we replace the double $p$-limit by the $p^2$-limit. We add a proof for the sake of completeness.

\begin{coro}
    A topological semigroup $X$ has an idempotent if and only if $f_x$ has a $p^2$-limit for some $p\in\w^*$.
\end{coro}

\begin{proof}
If $x\in X$ is idempotent, then $x$ is the $p^2$-limit of $f_x$ for every free ultrafilter $p\in\beta\w$. Hence assume that $e$ is the $p^2$-limit of $f_x$, we shall show that $e=e\cdot e$. Assume towards a contradiction that it is not the case and take $U\in \cN(e)$ such that $U^2\cap U=\emptyset$. Then $\{\{n,m\}:n<m\land x^{m-n}\in U\}\in p^2$. Let 
$$P_0=\{n\in\w:\{m\in\w:x^{m-n}\in U\}\in p\},$$
by assumption $P_0\in p$. Let $i=\min(P_0)$ and now define 
$$P_1=\{m\in P:x^{m-i}\in U\}\in p.$$
Take $j\in P_1$ and again define  
$$P_2=\{m\in P_1:x^{m-j}\in U\}\in p.$$
Finally, if $k\in P_2\subseteq P_1$, we get that $x^{k-i},x^{k-j},x^{j-i}\in U$ but $x^{k-i}=x^{k-j}\cdot x^{j-i}\in U^2$, which contradicts that $U$ and $U^2$ are disjoint.
\end{proof}

We will show that the previous Corollary really says more than Theorem \ref{banakhstheorem} by constructing a 2-countably compact space that is not doubly countably compact in Theorem \ref{2countablynotdoubly}.

Theorem \ref{banakhstheorem} is used in \cite{Banakhsemigroups} to obtain many properties that imply the existence of idempotents in topological semigroups and also to observe when topological semigroups became topological paragroups. Some of these properties have to do with powers of $X$ being countably compact, so the following question became relevant: 
\begin{problem}\label{problem}
Is the square of a doubly compact space countably compact?    
\end{problem}

We will consistently settle this question in the negative, by constructing a doubly compact space whose square is not countably compact under the assumption of the existence of a Ramsey ultrafilter. In addition, we will construct examples of a countably compact space that is not 2-countably compact and the already mentioned example of a 2-countably compact space that is not doubly countably compact.

Figure 2 shows the relations between 2 dimensional versions of countable compactness. Theorem \ref{doublywhosesquareisnot} will improve the classical result of Novak and Terasaka that there is a countably compact space whose square is not countably compact \cite{novakcartesianproduct}\cite{terasakacartesianproduct}. The price we pay for strengthening countable compactness to double countable compactness is that we need to assume the existence of Ramsey ultrafilters.

{\begin{figure}[ht]\label{Diagramindimension2}
    \centering
    \includegraphics[height=6cm]{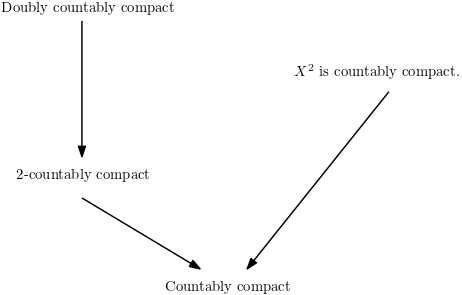}
    \caption{Relations between 2 dimensional versions of countable compactness}
\end{figure}}

\vspace{3mm}

Following \cite{Froliksumofuf}, we say that a family $\{q_n:n\in\w\}$ of ultrafilters on a countable set $N$ is \emph{discrete} if there is a partition $\{P_n:n\in\w\}$ of $N$ such that $P_n\in q_n$ for every $n\in\w$. Equivalently, $\{q_n:n\in\w\}$ is discrete if it is a discrete subset of $\beta N$, the Stone-\v{C}ech compactification of $N$ with the discrete topology. Given $p\in\w^*$, let also $\sum_p\{q_n:n\in\w\}$ denote the \emph{Frolik sum} of $\{q_n:n\in\w\}$ over $p$, given by
$$A\in\sum_p\{q_n:n\in\w\}\Leftrightarrow \{n\in\w: A\cap P_n\in q_n\}\in p.$$

We are mainly interested in ultrafilters on $\w\times\w$. Let $\pi_1$ and $\pi_2$ be the projections of $\w\times\w$ onto the first and second coordinates respectively. We can extend $\pi_1$ and $\pi_2$ to $\beta(\w\times\w)$ by declaring
$$\pi_i(p)=\{\pi_i[A]:A\in p\}$$
for $p$ an ultrafilter on $\w\times\w$. Thus $\pi_i(p)$ is an (possibly principal) ultrafilter on $\w$. Given two ultrafilters $p,q\in\w^*$ we define the \emph{Fubini product} $p\otimes q$ as the ultrafilter on $\w\times\w$ such that
$$A\in p\otimes q\Leftrightarrow \{n\in\w:\{m\in\w:(n,m)\in A\}\in q\}\in p.$$

Notice that Fubini products are essentially the same as Frolik sums where the sequence $q_i$ is such that $\{i\}\times\w\in q_i$ and $\pi_2(q_i)=q$ for every $i\in\w$. 

We say that an ultrafilter $q$ on $\w\times\w$ is a \emph{Fubini ultrafilter}, if there is an ultrafilter $p\in\w^*$ such that $q=p^2=p\otimes p$. The following lemmas can be deduced easily from the definitions and can be deduced from previous results in the literature (e.g., from Booth's article \cite{boothultrafilters}). However, we include them as we want to state them in a way that will be easier to manage on $\w\times\w$ and to get familiar with our notation and terminology.

\begin{lem}\label{uniquenessfroliksums}
    If $s_0\otimes r_0=q=s_1\otimes r_1$ then $s_0=s_1$ and $r_0=r_1$. 
\end{lem}

\begin{proof}
    Notice that if $q=s\otimes r$ then $s=\pi_1(q)$ is uniquely determined and so $s_0=s_1$. Let's call $s=s_0=s_1$. On the other hand, if $A\in r_0$, then $\w\times A\in s\times r_0=q=s\times r_1$, which implies $A\in r_1$. Hence $r_0\subseteq r_1$ and they are equal as both are ultrafilters.
\end{proof}

\begin{lem}\label{disjointnessfroliksums}
    Let $q\in\beta(\w\times\w)$, $r\in\w^*$ and $s=\pi_1(q)$. If $q\neq s\otimes r$ then there is an $A\in q$ such that $\pi_2(A\cap\{n\}\times\w)\notin r$ for every $n\in\w$. 
\end{lem}

\begin{proof}
Since $q\neq s\otimes r$ we can find a set $A'\in q$ such that 
$$\{n\in\w:\{m\in\w:(n,m)\in A'\}\in r\}\notin s$$
or equivalently
$$B=\{n\in\w:\{m\in\w:(n,m)\in A'\}\notin r\}\in s.$$
Thus $A=A'\cap(B\times\w)$ is as desired.
\end{proof}

\begin{lem}\label{froliksums1}
     If $\{i\}\times\w\in q_i$ for every $i\in\w$, then $\pi_1\big(\sum_s\{q_i:i\in\w\}\big)=s$. In particular $\pi_1(s\otimes r)=s$.
\end{lem}

\begin{proof}\label{froliksums2}
    It follows from the definitions.
\end{proof}

We now produce the first of our counterexamples.

\begin{thm}
    There is a space which is \emph{countably compact} but not $2$-\emph{countably compact}
\end{thm}
\begin{proof}

Let $X=\omega\times\omega$ and $\beta X$ be the Stone-Čech compactification of $\omega\times\omega$. As usual, we will denote by $X^*$ the remainder $\beta X$, that is, $X^*=\beta X\setminus X$.

We claim that  $Z=\beta X\setminus\{p\otimes p:p\in\omega^*\}$ is the desired space.
We first note that $Z$ is not $2$-\emph{countably compact} since the function $f:[\omega]^2\to Z$ where $f(\{n,m\})=(n,m)$, $n<m$ has no $p$-\emph{limit} for any $p\in\omega^*.$ 

Now we proceed to show that our space $Z$ is countably compact. To establish that, we need to show that every infinite set has an accumulation point, equivalently  every infinite sequence in $Z$ has a ${p}$-\emph{limit} for some $p\in\omega^*$.

\vspace{3mm}




Let $\langle{\bf{q}}_n:n\in\omega\rangle$ be a discrete sequence in $Z$. 
If there are infinitely many ${\bf{q}}_n$'s that are not principal, and thus ${\bf{q}}_n\in\w\times\w$, then we have an infinite subset $M\subseteq\w\times\w\cap \langle{\bf{q}}_n:n\in\omega\rangle$. We can further shrink $M$ and find $N\in[M]^\w$ such that either $N\subseteq \{i\}\times\w$ for some $i\in\w$ or $|N\cap \{i\}\times\w|\leq1$ for every $i\in\w$. Let $\bf{q}$ be any ultrafilter containing $N$, in any case $\bf{q}$ is a limit point of $\langle{\bf{q}}_n:n\in\w\rangle$ that is not a Fubini ultrafilter, hence ${\bf{q}}\in Z$. Then we can assume that ${\bf{q}}_n\in X^*$ for all $n\in\w$.
We investigate the possible cases:

\vspace{4mm}

\noindent {\bf CASE 1:} \underline{For infinitely many $n$, there is an $i_n$ such that $X_{i_{n}}=\{i_n\}\times\omega\in{\bf{q}}_n$.}

\vspace{4mm}

\noindent {\bf SUBCASE 1.1:} \underline{For infinitely many $n\in\w$, we have that $i_n=i$ for a fix $i\in\w$:}
In this subcase, $X_i\in {\bf{q}}_n$ for infinitely many $n\in\w$. Let $M=\{n\in\w:X_i\in {\bf{q}}_n\}$. As $\beta X$ is countably compact, we can find a limit point ${\bf{q}}\in\beta X$ of $\{{\bf{q}}_n:n\in M\}$. Notice that $X_i\in {\bf{q}}$ since otherwise $\{{\bf{r}}\in\beta X:(\w\times\w)\setminus X_i\}$ is an open neighborhood of ${\bf{q}}$ disjoint from $\{{\bf{q}}_n:n\in M\}$. Then it follows that ${\bf{q}}$ is not a Fubini ultrafilter as no Fubini ultrafilter can contain any $X_i$.
\vspace{4mm}

\noindent {\bf SUBCASE 1.2:} \underline{There is no such an $i\in\w$ as above.}
Then we can find an increasing sequence $\{i_n:n\in M\}$ for some infinite set $M\in[\w]^\w$ such that $X_{i_n}\in {\bf{q}}_n$ for every $n\in M$. Let $N=\{i_n:n\in M\}$, fix any ultrafilter $u\in\w^*\setminus\{\pi_2({\bf{q}}_n):n\in M\}$ such that $N\in u$ and define ${\bf{q}}=\sum_u\{{\bf{q}}_n:n\in M\}$, that is, $A\in {\bf{q}}$ if and only if
$$\{i_n:A\cap X_{i_n}\in {\bf{q}}_n\}\in u.$$
It is clear that ${\bf{q}}$ is a limit point of the sequence $\{{\bf{q}}_n:n\in M\}$.
As $\pi_1({\bf{q}})=u$, it suffices to show that ${\bf{q}}\neq u^2$ to conclude that ${\bf{q}}$ is not a Fubini ultrafilter. For this, note that since $u\neq \pi_2({\bf{q}}_n)$ for all $n\in M\in u$, we can find sets $A_n\in \pi_2({\bf{q}}_n)\setminus u$ such that
$$\bigcup_{n\in M}\{i_n\}\times A_n\in {\bf{q}}\setminus u^2.$$\\

\noindent {\bf CASE 2:} \underline{For all but finitely many $n\in\w$, no $X_i$ is an element of ${\bf{q}}_n$.}
\vskip 6pt
By trowing away finitely many ultrafilters in the sequence, we can assume that $X_i\not\in {\bf{q}}_n$ for every $n\in\w$. Let $\{Y_n:n\in\w\}$ be a partition of $\w\times\w$ such that $Y_n\in {\bf{q}}_n$ for every $n\in\w$. Notice that $\{X_n\cap Y_m:n,m\in\w\}$ is a partition of $\w\times\w$ such that  $Y_n\setminus\bigcup_{i<m}X_i\in {\bf{q}}_n$ for every $n,m\in\w$. 

    Consider the following condition:
    \begin{itemize}
        \item[($\star$)] There exists an ultrafilter ${\bf{r}}\in\w^*$ such that 
    $$M=\{n\in\w:{\bf{q}}_n=\pi_1({\bf{q}}_n)\otimes {\bf{r}})\}$$ 
    is infinite.
    \end{itemize}
    If $(\star)$ holds fix an ultrafilter ${\bf{r}}_0$ satisfying it. Let $M=M_0\sqcup M_1$ be a partition and pick $i<2$ such that $\bigcup_{n\in M_i}Y_n\notin {\bf{r}}_0\otimes {\bf{r}}_0$. Then fix an ultrafilter $u\in\w^*$ such that $N:=M_i\in u$.
    If condition $(\star)$ does not hold, let $u$ be any free ultrafilter on $\w$. Define ${\bf{q}}=\sum_u\{{\bf{q}}_n:n\in\w\}$. 

    It is clear that ${\bf{q}}$ is a limit of the sequence $\{{\bf{q}}_n:n\in\w\}$, so we are done if we show that ${\bf{q}}$ is not a Fubini ultrafilter.
    For this fix an ultrafilter ${\bf{r}}\in\w^*$. 
    If $(\star)$ holds and ${\bf{r}}={\bf{r}}_0$ it is clear that $\bigcup_{n\in N}Y_n\notin {\bf{r}}\otimes {\bf{r}}$ but it belongs to ${\bf{q}}$, so we can assume that either $(\star)$ does not hold or ${\bf{r}}\neq {\bf{r}}_0$.
    
    If $(\star)$ does not hold, then ${\bf{q}}_n\neq \pi_1({\bf{q}}_n)\otimes {\bf{r}}$ for all but finitely many $n\in\w$. 
    Even in the case where $(\star)$ does hold but ${\bf{r}}\neq {\bf{r}}_0$, by Lemma \ref{uniquenessfroliksums}, we can find a set $M\in u$ given by $(\star)$, such that ${\bf{q}}_n\neq \pi_1({\bf{q}}_n)\otimes {\bf{r}}$ for all $n\in M$. 
    In any case, by applying Lemma \ref{disjointnessfroliksums}, we can shrink $Y_n$ to a set $Q_n$ such that $Q_n\in {\bf{q}}_n$ and $\pi_2(Q_n\cap X_m)\notin {\bf{r}}$ for every $n\in M\in u$ and every $m\in \w$. As no $X_m\in {\bf{q}}_n$, we can further assume that $Q_n\cap \big(\bigcup_{i\leq n}X_i\big)=\emptyset$. Hence $D=\bigcup_{n\in M}Q_n\in {\bf{q}}$ and it remains to show that $C=(\w\times\w)\setminus D\in {\bf r}\otimes {\bf r}$.

    For $n\notin N$, let's take for simplicity $Q_n=\emptyset$.
    First of all, notice that 
    $$C=\bigcup_{n\in\w}X_n\setminus D=\bigcup_{n\in\w}\Big(X_n\setminus \bigcup_{i<n}Q_i\Big).$$
    Finally, $\pi_2\Big(X_n\setminus\bigcup_{i<n}Q_i\Big)\in {\bf{r}}$ for every $n\in\w$ since $\pi_2(Q_i\cap X_n)\notin {\bf{r}}$ for every $i,n\in\w$. Therefore $C\in {\bf{r}}\otimes {\bf{r}}$ as $\pi_2(C\cap X_n)=\pi_2\Big(X_n\setminus\bigcup_{i<n}Q_i\Big)\in{\bf{r}}$ for every $n\in\w$.
 \end{proof}

If we denote by $\M_2$ the class of Fubini ultrafilters on $\w\times\w$, the previous results says that $\beta(\w\times\w)\setminus\M_2$ is countably compact. It would be interesting to have a similar result for $n>2$. Or even better, to show that if $\M_{n+1}$ is defined similarly, then it is $n$-countably compact.
Formally, let $\M_n$ be the class of ultrafilters on $\w^n$ that are of the form $p^n=p\otimes\cdots\otimes p$. 

\begin{que}
    Let $n>2$. Is $\beta(\w^n)\setminus\M_n$ countably compact? Is it $m$-countably compact for every $m<n$?
\end{que}

In general, it would be nice to distinguish the classes of $\alpha$-countably compact spaces.

\begin{que}
    Are there $\alpha$-countably compact spaces that are not $\beta$-countably compact for $1<\alpha<\beta<\omega_1$?
\end{que}

\section{Ramsey ultrafilters}

In this section, we solve Problem \ref{problem} by making use of Ramsey ultrafilters. 
Recall that an ultrafilter $p\in\w^*$ is {\it Ramsey} if for every partition $\{P_n:n\in\w\}$ either $P_i\in p$ for some $i\in\w$ or there exists $U\in p$ such that $|U\cap P_i|=1$ for every $i\in\w$.
It is known that under \textsf{CH} there are $2^\c$ many Ramsey ultrafilters \cite{comfortultrafilters}. 
On the other hand, if there is a Ramsey ultrafilter, then it can be easily shown that there are $\c$-many: Just take a Ramsey ultrafilter $\cU$ an almost disjoint family $\{A_\alpha:\alpha<\c\}$ and for every $\alpha<\c$ a bijection $f_\alpha$ between $\w$ and $A_\alpha$. 
Hence $\{f_\alpha(\cU):\alpha<\c\}$ is a family of $\c$ pairwise distinct Ramsey ultrafilters. It was also shown by Baumgartner and Laver that it is consistent that there are exactly $\c$-many Ramsey ultrafilters \cite{baumgartneriterated}.

We recall also that an ultrafilter $u\in\beta\w$ is a $P$\emph{-point}, if every sequence $\{U_n:n\in\w\}\subseteq u$, has a pseudointersection $P\in u$. Another relevant class of ultrafilters in our constructions will be that of discrete ultrafilters.

\begin{dfn}\cite{banakhblass}
    We say that an ultrafilter $u\in\beta\w$ is $Y$-\emph{discrete} if for every one-to-one sequence $\{x_n:n\in\w\}\subseteq Y$, there exists $A\in u$ such that $\{x_n:n\in A\}$ is a discrete subset of $Y$.
\end{dfn}

If we let $Y=\bR$, we get the classical notion of a discrete ultrafilter due to Baumgartner \cite{baumgartnerultrafilters}. 
It is easy to see that every $P$-point is discrete, and moreover, every discrete ultrafilter is $\beta\w$-discrete (see \cite{banakhblass}).

It is also known that if there is a $P$-point, then there is a discrete ultrafilter that is not a $P$-point. On the other hand, it is not known if the existence of discrete ultrafilters implies the existence of $P$-points, and neither is it known whether $\beta\w$-discrete ultrafilters exist in {\sf ZFC}.

\begin{propo}\label{ramseydiscrete}
    Given a Ramsey ultrafilter $r$ and a sequence $\{x_n:n\in\w\}\subseteq\beta\w$, there exists $M\in r$ such that $\{x_n:n\in M\}$ is either constant or discrete.
\end{propo}

\begin{proof}
    It is well known that for a Ramsey ultrafilter, we can find such a subsequence that is either constant or one-to-one. As every $P$-point is $\w^*$-discrete the result follows.
\end{proof}



\begin{nota}
    We say that a double sequence $\{x_n^m:n,m\in M\}\subseteq X$ is discrete, where $M\in[\w]^\w$, if $f:[M]^2\to X$ given by $f(\{n,m\})=x_n^m$, for $n<m$, is one-to-one and its range is discrete
\end{nota}

It is easy to see that given a discrete double sequence on $\beta\w$, the $p^2$-limit and the $q^2$-limit are distinct for $p\neq q$. Thus discrete double sequences have $2^\c$-many $p^2$-limit points. 
In order to construct a space which is doubly countably compact but whose square is not countably compact, we need to show that any double sequence is either somewhere constant or it has many disjoint $p$-limit sequences. The first step, is to show that for discrete double sequences, we also have $2^\c$-many pairwise disjoint $p$-limit sequences.
For the next lemma, given $M\in[\w]^\w$, we remind that $\beta M$ coincides with the subset of ultrafilters $p\in\beta\w$ such that $M\in p$. 

\begin{lem}\label{discretesquare}
    Let $M\in[\w]^\w$ and let $\{x_n^m:n,m\in M\}\subseteq\beta\w$ be a discrete double sequence. Then for every set $A\subseteq\beta\w$ of size less than $2^\c$, there exists a $p$-limit sequence $(x,x_n)_{n\in M}$ for $\{x_n^m:n,m\in M\}$ which is disjoint from $A$. 
\end{lem}

\begin{proof}
    Let $x_n^p\in\beta\w$ be the $p$-limit of $\{x_n^m:m\in M\}$ for each $n\in M$ and $p\in\beta M$. In addition, for every $n\in M$, let $F(n)=\{p\in\beta M:x_n^p\in A\}$. As $\{x_n^m:m\in M\}$ is discrete, we have that $|F(n)|\leq|A|<2^\c$ and as $2^\c$ has uncountable cofinality, $|\bigcup_{n\in M}F(n)|<2^\c$ as well. Let now $x^p=p^2-\lim\{x_n^m:n,m\in M\}$ and define $G=\{p\in\beta M:x^p\in A\}$. As $p^2$-limits are unique, we have that $|G|\leq A<2^\c$. Finally, we see that for every $p\in\beta M\setminus (F\cup G)$, the $p$-limit sequence $(x^p,x_n^p)_{n\in M}$ is disjoint from $A$.
\end{proof}

We now proceed to prove our main Lemma.

\begin{lem}\label{mainlemma}
    Assume that there are $\kappa$-many Ramsey ultrafilters. If $A\subseteq\beta\w$ and $|A|<\kappa$, then $\beta\w\setminus A$ is doubly countably compact.
\end{lem}

\begin{proof}
    Let $\{x_n^m:n<m<\w\}\subseteq\beta\w\setminus A$. The goal is to find a $p$-limit sequence for this sequence which is disjoint from $A$. We can find a decreasing sequence $\{M_n:n\in\w\}$ such that $\{x_n^m:m\in M_n\setminus(n+1)\}$ is either constant or discrete. As $x_n^m$ is defined only when $m>n$, we can find a pseudointersection $M$ such that $\{x_n^m:m\in M\setminus(n+1)\}$ is either constant or discrete for every $n\in M$.

    \vspace{3mm}

    \noindent{}{\bf CASE 1: }$\{x_n^m:m\in M\setminus(n+1)\}$ is constant for infinitely many $n\in M$.

    Let $M_0$ be the set of these $n$ and let $x_n$ be the constant value that takes on the sequence $\{x_n^m:m\in M_0\setminus(n+1)\}$. We remark that $x_n\notin A$ for every $n\in M_0$ and that $x_n=p$-$\lim\{x_n^m:m\in M_0\setminus(n+1)\}$ for every $p\in\w^*$ containing $M_0$. Take $M_1\in[M_0]^\w$ such that now $\{x_n:n\in M_0\}$ is either constant or discrete.

    \vspace{3mm}

    \noindent {\bf Subcase 1.1: }$\{x_n:n\in M_0\}$ is constant.

    Assume that the value of this sequence is $x$, then $(x,x_n)_{n\in M_1}$ is the $p$-limit sequence of $\{x_n^m:m\in M_1\setminus(n+1)\}$ for every $p\in\w^*$ containing $M_1$. As $x=x_n\notin A$ we have constructed a $p$-limit sequence disjoint from $A$.

    \vspace{3mm}

    \noindent {\bf Subcase 1.2: }$\{x_n^m:m\in M_1\setminus(n+1)\}$ is discrete.

    As $|\overline{\{x_n:n\in M_1\}}|=2^\c$, we can find $p\in\w^*$ with $M_1\in p$ and $x\notin A$ such that $x=p$-$\lim\{x_n:n\in M_1\}$. Then $x$ is the double $p$-limit of $\{x_n^m:m\in M_1\setminus(n+1)\}$ through $\{x_n:n\in M_1\}$ and this $p$-limit sequence is disjoint from $A$.

    \vspace{3mm}

    \noindent {\bf CASE 2: } Without loss of generality, $\{x_n^m:m\in M\setminus(n+1)\}$ is discrete for every $n\in M$.

    For each $n\in M$ and each ultrafilter $p\in\w^*$ containing $M$, let $x_n^p$ be the $p$-limit of $\{x_n^m:m\in M\setminus(n+1)\}$. If $p$ is Ramsey, find $M_p\in[M]^\w$ such that $\{x_n^p:n\in M_p\}$ is either constant or discrete. 

    \vspace{3mm}
    
    \noindent {\bf Subcase 2.1: }There is a Ramsey ultrafilter $p$ such that $\{x_n^p:n\in M_p\}$ is discrete.

    Let $\{Q_n:n\in M_p\}$ be a partition of $\w$ such that $Q_n\in x_n$ for every $n\in M_p$. Then we have that $A_n=\{m\in\w:Q_n\in x_n^m\}\in p$ for each $n\in M_p$. 
    Define $m_n=\min\bigcap_{i\leq n}A_i$ and let $M_1=\{m_n:n\in\w\}$. Without loss of generality we can assume that $M_1\subseteq M_p\subseteq M$, so $\{x_n^m:m\in M_1\setminus(n+1)\}$ is discrete for every $n\in M_1$. 
    Moreover, we have that $Q_n\in x_n^m$ for each $n\in M_1$ and each $m>n$. 
    Combining these two observations, we can find further partitions $\{Q_n^m:m\in M_1\setminus(n+1)\}$ for $n\in M_1$ such that $Q_n^m\in x_n^m$. Therefore the partition $\{Q_n^m:n,m\in M_1\land m>n\}$ witnesses that the sequence $\{x_n^m:n,m\in M_1\land m>n\}$ is discrete and we can find a $p$-limit sequence disjoint from $A$ by Lemma \ref{discretesquare}.

    \vspace{3mm}
    \noindent {\bf Subcase 2.2: }$\{x_n^p:n\in M_p\}$ is constant for every Ramsey ultrafilter $p$.
    Let $x^p$ be this constant value for any $p$ Ramsey. Note that if there is $x\in\w^*$ such that 
    $$|\{p\in\w^*:p\textnormal{ is Ramsey and }x^p=x\}|\geq\w_1,$$
    then for some $n\in\w$ and for two distinct ultrafilters Ramsey $p,q\in\w^*$ we obtain that $n\in M_p\cap M_q$ and thus $x_n^p=x=x_n^q$. But this is a contradiction as
    \begin{itemize}
        \item $x_n^p=p$-$\lim\{x_n^m:m\in M\setminus(n+1)\}$
        \item $x_n^q=q$-$\lim\{x_n^m:m\in M\setminus(n+1)\}$
    \end{itemize}
    and the sequence $\{x_n^m:m\in M\setminus(n+1)\}$ is discrete.
    
    Thus, every $x\in\beta\w$ equals $x^p$ for at most countably many $p$ Ramsey ultrafilters. Let $p$ be any Ramsey ultrafilter such that $x^p\notin A$ (this is possible as there are only $|A|\cdot\w<\kappa$ forbidden ultrafilters). Thus the $p$-limit sequence $(x^p,x_n^p)_{n\in M_p}$ is constant with value $x^p\notin A$ and we are done.
\end{proof}

As mentioned above, the existence of a Ramsey ultrafilter implies the existence of at least $\c$-many of them. Hence we have the following two corollaries in the extreme cases for the number of existing Ramsey ultrafilters.

\begin{coro}
    If there is a Ramsey ultrafilter and $A\subseteq\beta\w$ has size less than $\c$, then $\beta\w\setminus A$ is doubly countably compact.
\end{coro}

\begin{coro}
    If there are $2^\c$-many Ramsey ultrafilters, $A\subseteq\beta\w$ and $|A|<|\beta\w|$, then $\beta\w\setminus A$ is doubly countably compact.
\end{coro}

Our next result can be seen as a generalization of the classical result of Novak \cite{novakcartesianproduct} and Terasaka \cite{terasakacartesianproduct} that shows that a countably compact space need not to have countably compact square. The price we pay for strengthening countable compactness to double countable compactness, is that we only have a consistent result by assuming the existence of a Ramsey ultrafilter.

\begin{thm}\label{doublywhosesquareisnot}
    Assume there is a Ramsey ultrafilter. Then there is a doubly countably compact space whose square is not countably compact.
\end{thm}

\begin{proof}
We construct doubly countably compact spaces $X,Y\subseteq\beta\w$ such that $X\cap Y=\w$. Hence $\{(n,n)\in X\times Y:n\in\w\}$ has no accumulation points in $X\times Y$, as any accumulation point $z=(x,y)$ would be the $p$-limit of $(n,n)_{n\in\w}$ for some $p\in\w^*$, but this implies $x=p=p-\lim\w=y$.
In other words, $(x,y)\in X\times Y$ and $x=y=p\in X\cap Y\setminus \w$, which is a contradiction.

To construct $X$ and $Y$ we perform a recursion of length $\c$. Let $\{f_\alpha:\alpha<\c\}$ enumerate all functions $f:[\w]^2\to\w\cup\w\times[\w,\c)$. We start by defining $X_\w=\w=Y_\w$. In general for an infinite ordinal $\alpha<\c$, let
$$X_\alpha=\bigcup_{\beta<\alpha}X_\beta$$
and
$$Y_\alpha=\bigcup_{\beta<\alpha}Y_\beta$$
for $\alpha\in\c$ limit.

If $\alpha=\beta+1$, then 
$$X_\alpha=X_\beta\cup\{x_\beta^n:n\in\w\}$$ 
for some $\{x_\beta^n:n\in\w\}\subseteq\beta\w$ and similarly $$Y_\alpha=Y_\beta\cup\{y_\beta^n:n\in\w\}.$$ 
We will make sure that $X_\alpha\cap Y_\alpha=\w$ for $\alpha<\c$ along the construction.

Assume we have constructed $X_\beta$ and $Y_\beta$ for $\alpha=\beta+1<\c$. Consider the functions $g_\beta:[\w]^2\to X_\beta$ and $h_\beta:[\w]^2\to Y_\beta$ where 
\[
  g_\beta(\{n,m\})=
  \begin{cases}
    k & \text{if $f_\beta(\{n,m\})=k$} \\
    x_\alpha^k & \text{if $f_\beta(\{n,m\})=(k,\gamma)$}
  \end{cases}
\]
and $h_\beta$ is defined similarly with $y_\gamma^k$ instead of $x_\gamma^k$. Let $Y_\beta'=Y_\beta\setminus\w$. As $|Y_\beta'|\leq\w\cdot\beta<\c$, we can apply Lemma \ref{mainlemma} and find a $p$-limit sequence $(x,x_n)_{n\in\w}$ for $g_\beta$ which is disjoint from $Y_\beta'$. Let $x_\beta^0=x$ and $x_\beta^n=x_n$ for every $n>0$. Hence $X_\alpha$ as defined above contains a doubly $p$-limit for $g_\beta$. We define similarly $\{y_\beta^n:n\in\w\}$. 

Finally let $X=\bigcup_{\alpha<\c}X_\alpha$ and $Y=\bigcup_{\alpha<\c}Y_\alpha$. It is clear from our inductive assumptions that $X\cap Y=\w$. Moreover, $X$ and $Y$ are doubly countably compact as any double sequence in each of them appears as $g_\beta$ or $h_\beta$ at some stage $\beta<\c$. Hence we have two doubly countably compact spaces whose product is not countably compact. By taking their disjoint sum, we get the desired doubly countably compact space whose square is not countably compact.
\end{proof}

By assuming the stronger hypothesis that there are more than $\c$ selective ultrafilters we can greatly simplify the previous proof as follows: Take an elementary submodel $M$ of size $\c$ of $H(\theta)$ for $\theta$ large enough. Thus $X=\beta\w\cap M$ is doubly countably compact and so is $Y=\beta\w\setminus A'$, where $A'=A\setminus\w$. Therefore, following the first paragraph in the proof of the previous Theorem we get the conclusion.

By closely looking at the proof, it is clear that we can moreover construct countably many, doubly countably compact subspaces of $\beta\w$ whose intersection of every pair is $\w$. Similarly, by assuming the existence of more than $\c$ Ramsey ultrafilters and by bookkeeping argument, we can construct $\c$-many. So we have the following corollary.

\begin{coro}\label{disjointfamilies}
    If there is a Ramsey ultrafilter, then there is a countable family $\cD$ of doubly countably compact subspaces of $\beta\w$ such that $A\cap B=\w$ for every distinct $A,B\in \cD$. Assuming there are more than $\c$-many Ramsey ultrafilters, we can moreover assume that $|\cD|=\c$. 
\end{coro}

We now show that 2-countable compactness is consistently weaker than double countable compactness.

\begin{thm}\label{2countablynotdoubly}
    Assume there is a Ramsey ultrafilter. Then there is a 2-countably compact space that is not doubly countably compact.
\end{thm}

\begin{proof}
For ease of notation, we work on $\beta(\w\times\w)$ again. Let $B_n=\{n\}\times\w$ for every $n\in\w$. 
Thus $\overline{B}_n\cap\overline{B}_m=\emptyset$ whenever $n\neq m$. 

By Corollary \ref{disjointfamilies}, there is a family $\{A_n\subseteq\beta\w:n\in\w\}$, such that $A_n\cap A_m=\w$ for $n\neq m$ and each $A_n$ is doubly countably compact. 
Let $C_n$ be the homeomorphic copy of $A_n$ in $\overline{B_n}$ given by the natural enumeration of $B_n$. 
Thus $\{C_n:n\in\w\}$ is a countable family of disjoint doubly countably compact subspaces of $\beta(\w\times\w)$ and moreover, if $x_n^p$ is the $p$-limit of $B_n=\{n\}\times\w$ for $p\in\w^*$ and $n\in\w$, then $x_n\notin C_n$ for all but at most one $n\in\w$. 

Let $C_\w=\{p\in\beta(\w\times\w):\forall k\in\w\ (B_k\notin p)\}$ and define $X=\bigcup_{n\in\w+1}C_n$. By the observations made on $\{C_n:n\in\w\}$ it is clear that the function $G:[\w]^2\to X$ given by $G(\{n,m\})=(n,m)$ for $n<m$ does not have a double $p$-limit in $X$, so $X$ is not doubly countably compact.

It remains to show that $X$ is 2-countably compact. Let $f:[\w]^2\to X$ be arbitrary. 
By Ramsey's theorem, find $M_{0}\in[\w]^\w$ such that $f''[M_0]^2$ is either contained or disjoint from $C_0$. If the former case holds, let $M=M_0$, otherwise, we can repeat the argument with $C_1$ and define a set $M_1\in[M_0]^\w$. Continuing in this way, if for some $i\in\w$ we have $f''[M_i]^2\subseteq C_i$, let $M=M_i$, otherwise, we have a decreasing sequence $\{M_i:i<\w\}$ and thus we can define a pseudointersection $M$ such that $f''[M\setminus k]^2\subseteq\bigcup_{n\in(\w+1)\setminus k}C_n$.

\vspace{3mm}

If $f''[M]^2\subseteq C_j\rest\beta$ for some $j\in\w$, fix this $j$. As $C_j$ is doubly countably compact (so 2-countably compact), we can find a $p^2$-limit for $f\rest[M]^2$ and hence a $q^2$-limit for $f$, for some $q$ containing $M$.

Otherwise, $M$ satisfies that $f''[M\setminus k]^2\subseteq\bigcup_{n\in(\w+1)\setminus k}C_n$. Since $\beta(\w\times\w)$ is 2-countably compact, there is a $p\in\w^*$ such that $x=p^2-\lim f\in\beta(\w\times\w)$. 
We claim that $x\in C_\w$. To see this, notice that if $B_k\in x$ for some $k\in\w$, then $U=\{z\in\beta(\w\times\w):B_k\in z\}$ is an open neighborhood of $x$ and 
$$\{\{n,m\}:f(\{n,m\})\in U\}\subseteq\{\{n,m\}:\min(\{n,m\})<k\}\in\sf{FIN}^2,$$
which is impossible as $p$ is a free ultrafilter and $x$ is the $p^2$-limit of $f$. Therefore $x\in C_\w\subseteq X$.
\end{proof}

Our proof of Lemma \ref{mainlemma} was carried out in {\sf ZFC} in its first case. For the second case, we need to assume the existence of Ramsey ultrafilters in order to find a subsequence which is either constant or discrete. 
If it was the case that for every $P$-point $p$ and every sequence in $\w^*$ we can find a subsequence which is either constant or one-to-one indexed by a set in $p$, then this can be easily accomplished, as every $P$-point is a discrete ultrafilter. However, it is easy to see that it is not the case if $p$ is a non-Ramsey $P$-point.

Perhaps, a more sophisticated construction would lead to the same result under weaker assumptions, or even, in {\sf ZFC}. We will introduce the class of $(\beta\w,2)$-discrete ultrafilters, which is a reasonable candidate for the class of ultrafilters that would be useful for this task.

\begin{dfn}
    Given a topological space $Y$ and $n\in\w$, we say that $u\in\beta\w$ is \emph{$(Y,n)$-discrete}, if for every one-to-one map $f:[\w]^n\to \beta\w$, there is $M\in u$ such that $f''[M]^n$ is discrete.
\end{dfn}

The previous definition can be extended to countable ordinal $\alpha$ by considering the quantification: ``for all barrier $\cB$ of rank less than or equal to $\alpha$ and for all $f:\cB\to Y$''.

We are mainly interested in $(\beta\w,2)$-discrete ultrafilters, so we are tempted to simply call hem 2-discrete, however, for historical reasons, we will save the term $n$\emph{-discrete ultrafilter} in order to refer to $(\bR,n)$-discrete ultrafilters.

It follows easily from Proposition \ref{ramseydiscrete} and the argument in Subcase 2.1 of Lemma \ref{mainlemma} that Ramsey ultrafilters are $(\beta\w,2)$-discrete and clearly $(\beta\w,2)$-discrete ultrafilters are $\beta\w$ discrete. Hence $(\beta\w,2)$-discrete ultrafilters sit between Ramsey and $\beta\w$-discrete ones just as $P$-points do. Figure 3 visualizes the relations just mentioned. The following questions are then natural.

\begin{que}
    Is every $P$-point $(\beta\w,2)$-discrete? What about the converse?
\end{que}

\begin{que}
    Do $(\beta\w,2)$-discrete ultrafilters coincide with Ramsey ultrafilters?
\end{que}

It would also be nice to show that Ramsey ultrafilters are 2 discrete.

\begin{que}
    Is every Ramsey ultrafilter 2-discrete? Is it $n$-discrete for every $n\in\w$?
\end{que}

{\begin{figure}[ht]\label{Diagramultrafilters}
    \centering
    \includegraphics[height=5cm]{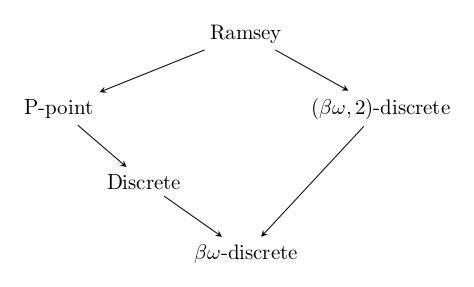}
    \caption{Special kinds of ultrafilters}
\end{figure}}

Of course our main goal is to produce the examples in Theorem \ref{doublywhosesquareisnot} and Theorem \ref{2countablynotdoubly} under weaker assumptions, so we ask:

\begin{que}
    Is there a doubly countably compact space whose square is not countably compact in {\sf ZFC}? Under the existence of $P$-points or $(\beta\w,2)$-discrete ultrafilters?
\end{que}

\begin{que}
    Is there a 2-countably compact but not doubly countably compact space  in {\sf ZFC}? Under the existence of $P$-points or $(\beta\w,2)$-discrete ultrafilters?
\end{que}

Such examples could be easily deduced from a positive answer to the following:

\begin{que} Assume $X\subseteq \omega^*$ has size ${\mathfrak c}$. Is $\beta\omega\setminus X$  2-countably compact?
\end{que}

Even a positive answer assuming $|X|<{\mathfrak c}$ would suffice. Also, we don't know if $\beta\omega\setminus X$ would be $\cB$-countably compact whenever $X$ is of size $\leq {\mathfrak c}$.

In the other direction, we do not know if a space having countably compact square implies that it satisfies a stronger version of countable compactness, like the 2 dimensional versions considered in this paper.

\begin{que}
    If $X^2$ is countably compact, does it follow that $X$ is doubly countably compact or at least 2 countably compact?
\end{que}

Finally, a positive answer to the following question would give us the above examples in ZFC:

\begin{que} If $f:[\omega]^2\rightarrow \beta\omega$ is given. Is there always an infinite $M$ such that $f[M]^2$ either is constant or has $2^c$ many $p^2$ limit points?
\end{que}

\subsection*{Acknowledgments}
The first author thanks Serhii Bardyla for pointing out Theorem \ref{banakhstheorem}.

\bibliography{convergence.bib}{}
\bibliographystyle{plain}

\end{document}